\documentclass[reqno]{amsart}
\usepackage[centertags]{amsmath}
\usepackage{graphicx}
\usepackage{amsfonts}
\usepackage{amssymb}
\usepackage[top=1in, bottom=1in, left=1in, right=1in]{geometry}

\usepackage[pdfborder={0 0 0},
   pdftitle={A FINITE SUBDIVISION RULE FOR THE N-DIMENSIONAL TORUS},
  pdfauthor={BRIAN RUSHTON}, pdftex, bookmarks=true, bookmarksnumbered = true]{hyperref}

\theoremstyle{definition}
\newtheorem{thm}{Theorem}
\theoremstyle{definition}

\theoremstyle{definition}
\newtheorem{defi}{Definition}
\theoremstyle{remark}

\theoremstyle{definition}

\begin{document}
\date{\today}

\title{A finite subdivision rule for the n-dimensional torus.}%

\author{Brian Rushton}

\maketitle
\pdfbookmark[1]{A FINITE SUBDIVISION RULE FOR THE N-DIMENSIONAL TORUS}{user-title-page}

\begin{center}



\end{center}

\begin{abstract}
Cannon, Floyd, and Parry have studied subdivisions of the 2-sphere extensively, especially those corresponding to 3-manifolds, in an attempt to prove Cannon's conjecture. There has been a recent interest in generalizing some of their tools, such as extremal length, to higher dimensions. We define finite subdivision rules of dimension $n$, and find an $n-1$-dimensional finite subdivision rule for the $n$-dimensional torus, using a well-known simplicial decomposition of the hypercube.
\end{abstract}

\section{Introduction}\label{Introduction}

In its most general form, a subdivision rule it is an algorithm for taking a surface with a labeled finite covering by compact sets and recursively refining the elements of the covering into smaller labeled compact sets.

Cannon and Swenson have shown \cite{hyperbolic} that Gromov hyperbolic groups with a 2-sphere at infinity give rise in a natural way to a subdivision rule on the sphere. In this setting, the subdivision rules are allowed to act on coverings of the 2-sphere by overlapping compact sets. A \emph{finite} subdivision rule is a simpler form of subdivision rule that other subdivision rules can often be reduced to, in which the coverings have to be tilings. Such subdivision rules have been studied extensively by Cannon, Floyd, and Parry in an attempt to solve the following conjecture \cite{combinatorial}:

\textbf{Conjecture:} All Gromov hyperbolic groups with a 2-sphere at infinity act co-compactly and properly discontinuously on $\mathbb{H}^3$ by isometries.

Finite subdivision rules are only one tool in studying this conjecture. Recently, several researchers have been expanding some of the other tools of Cannon, Floyd and Parry to 3-dimensions, such as extremal length and sphere-packings (see \cite{Spherepackings}, \cite{Spherelackings}, and \cite{saar}).

\begin{figure}
\begin{center}
\scalebox{1.4}{\includegraphics{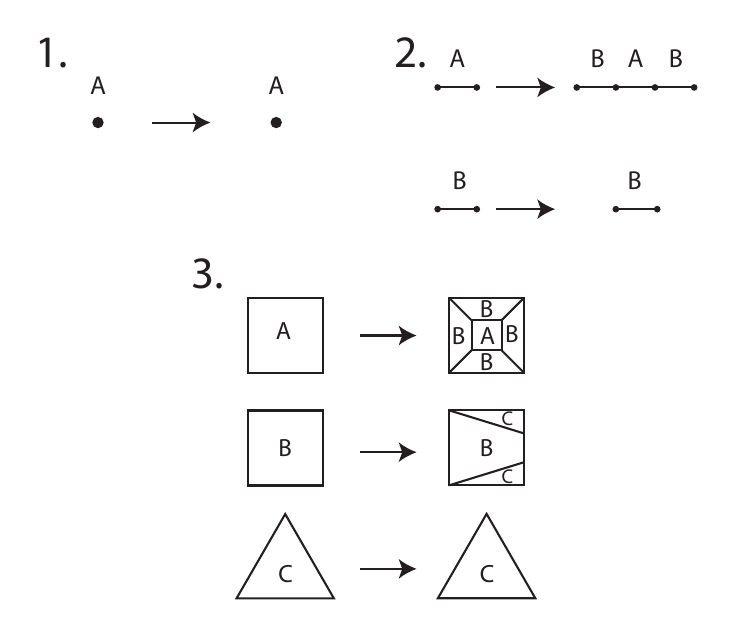}} \caption[The subdivision rules for the 1-,2- and 3-torus.]{1. The subdivision rule for the 1-torus, i.e. the
circle. 2. The subdivision rule for the 2-torus, i.e. the standard torus. 3. The subdivision rule for the 3-torus.} \label{SmallToriSubs1}
\end{center}
\end{figure}

\begin{figure}
\begin{center}
\scalebox{.8}{\includegraphics{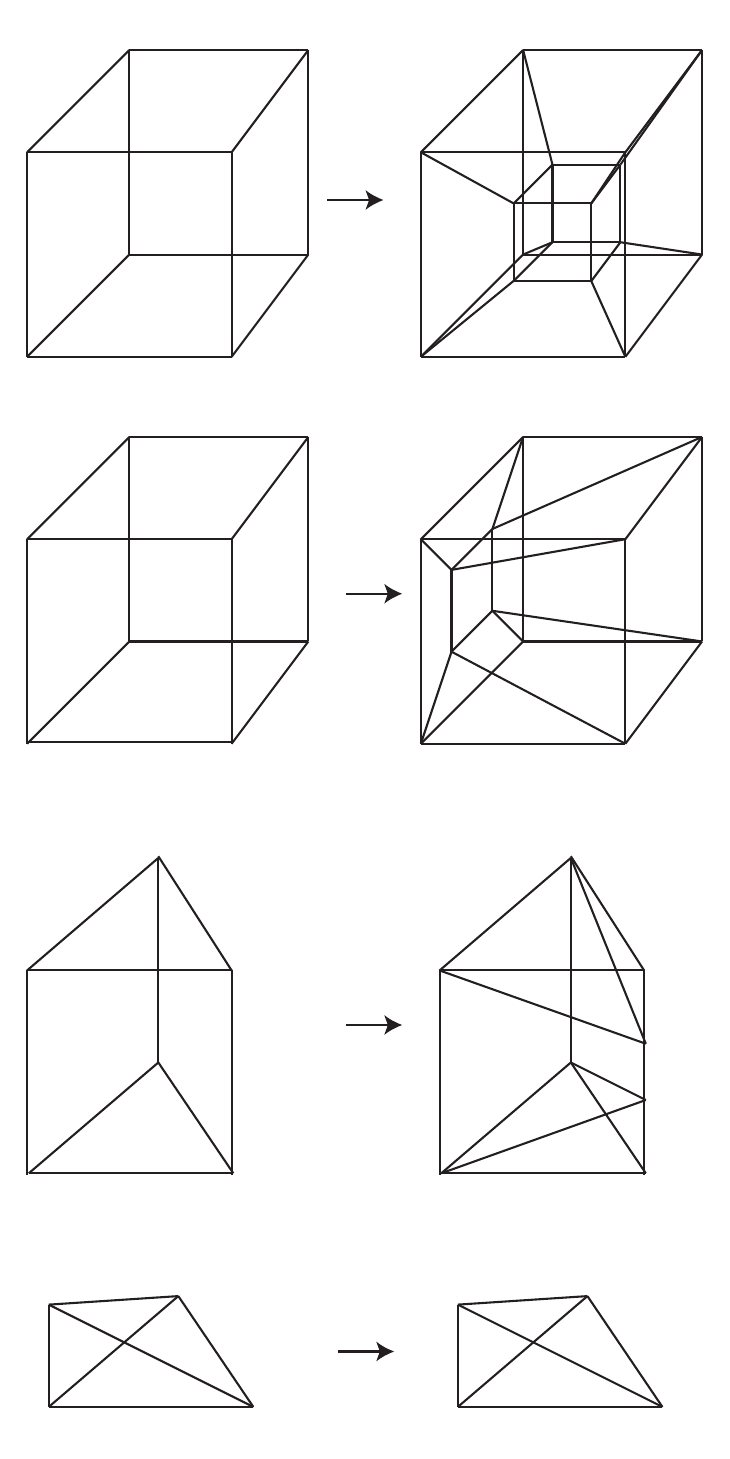}} \caption[The tile types for the four-torus.]{The tile types for the four-torus. The specifics about orientation and tile types will be discussed in more detail in the proof of Theorem \ref{CubeTheorem}.} \label{FourTorusSubs1}
\end{center}
\end{figure}

In this same spirit, the author has investigated subdivision rules of the 3-sphere arising from boundaries of 4-manifolds. The least complicated non-simply connected 4-manifold is the 4-torus $S^1 \times S^1 \times S^1 \times S^1$. With some effort, the method of Chapter 2 of \cite{myself2} for finding a finite subdivision rule for the 3-torus (pictured in part 3. of Figure \ref{SmallToriSubs1}) generalizes to give a subdivision rule for the 4-torus (shown in Figure \ref{FourTorusSubs1}). Notice that the tile types in Figure \ref{FourTorusSubs1} are (as CW-complexes) a a cube, another cube, a triangular prims, and a tetrahedron. We can state this more simply by letting $\Delta^n$ represent a simplex of dimension $n$, and $I^m$ a hypercube of dimension $m$. With this notation, these tile types are $I^3, I^2 \times \Delta, I \times \Delta^2,$ and $\Delta^3$. A quick glance at the finite subdivision rule for the 3-torus in Figure \ref{SmallToriSubs1} shows that the tile types here are $I^2, I \times \Delta,$ and $\Delta^2$. These two examples follow a definite pattern.

The main purpose of this paper is to show that this pattern continues for all $n$. In Theorem \ref{CubeTheorem} on page \pageref{CubeTheorem}, we construct an explicit $n$-dimensional finite subdivision rule (which we define below) for the $(n+1)$-torus in which the tile types have the form $I^k \times \Delta^{n-k}$, with one tile type for each $k$ from $0$ to $n$. Parts 1. and 2. of Figure \ref{SmallToriSubs1} show that this pattern holds for $n=1$ and $n=2$, as well.

\section{Formal Definition of a Subdivision Rule}

At this point, it will be helpful to give a concrete definition of subdivision rule. We first recall Cannon, Floyd and Parry's definition of a finite subdivision rule, taken from \cite{subdivision}.

\begin{defi} A \textbf{finite subdivision rule} $R$ consists of the following.
\begin{enumerate}
\item A finite 2-dimensional CW complex $S_R$, called the \textbf{subdivision complex}, with a fixed cell structure such that $S_R$ is the union of its closed 2-cells. We assume that for each closed 2-cell $\tilde{s}$ of $S_R$ there is a CW structure $s$ on a closed 2-disk such that $s$ has at least three vertices, the vertices and edges of $s$ are contained in $\partial s$, and the characteristic map $\psi_s:s\rightarrow S_R$ which maps onto $\tilde{s}$ restricts to a homeomorphism onto each open cell.
\item A finite two dimensional CW complex $R(S_R)$, which is a subdivision of $S_R$.
\item A continuous cellular map $\phi_R:R(S_R)\rightarrow S_R$ called the \textbf{subdivision map}, whose restriction to every open cell is a homeomorphism.
\end{enumerate}
\end{defi}

Each CW complex $s$ in the definition above (with its given characteristic map $\psi_s$) is called a \textbf{tile type}.

As the final part of the definition, they show how finite subdivision rules can act on surfaces (and 2-complexes in general). An $R$-complex for a subdivision rule $R$ is a 2-dimensional CW complex $X$ which is the union of its closed 2-cells, together with a continuous cellular map $f:X\rightarrow S_R$ whose restriction to each open cell is a homeomorphism. We can subdivide $X$ into a complex $R(X)$ by requiring that the induced map $f:R(X)\rightarrow R(S_R)$ restricts to a homeomorphism onto each open cell. $R(X)$ is again an $R$-complex with map $\phi_R \circ f:R(X)\rightarrow S_R$. By repeating this process, we obtain a sequence of subdivided $R$-complexes $R^n(X)$ with maps $\phi_R^n\circ f:R^n(X)\rightarrow S_R$. All of the preceding definitions were adapted from \cite{subdivision}, which contains several examples. While in theory, a subdivision rule is represented by a CW-complex, most rules in practice are described by diagrams of the sort shown in Figure \ref{SmallToriSubs1}, part 3.

In this paper, we find subdivision rules for the $n+1$-dimensional torus which subdivide the $n$-sphere. We define a subdivision rule in higher dimensions in a way analogous to subdivision rules in dimension 2. A \textbf{finite subdivision rule $R$ of dimension $n$} consists of:
\begin{enumerate}
\item A finite $n$-dimensional CW complex $S_R$, called the \textbf{subdivision complex}, with a fixed cell structure such that $S_R$ is the union of its closed $n$-cells. We assume that for every closed $n$-cell $\tilde(s)$ of $S_R$ there is a CW structure $s$ on a closed $n$-disk such that any two subcells that intersect do so in a single cell of lower dimension, the subcells of $s$ are contained in $\partial s$, and the characteristic map $\psi_s:s\rightarrow S_R$ which maps onto $\tilde{s}$ restricts to a homeomorphism onto each open cell.
\item A finite $n$-dimensional subdivision $R(S_R)$ of $S(R)$.
\item A \textbf{subdivision map} $\phi_R: R(S_R)\rightarrow S_R$, which is a continuous cellular map that restricts to a homeomorphism on each open cell.
\end{enumerate}

Each CW complex $s$ in the definition above (with its appropriate characteristic map) is called a \textbf{tile type} of $S$. All other portions of the definition (such as $R$-complexes) generalize in the natural way. As for traditional finite subdivision rules, we will often describe an $n$-dimensional finite subdivision rule by the subdivision of every tile type, instead of by constructing an explicit complex.

Given a finite subdivision rule $R$ of dimension $n$, an $R$-complex consists of a $n$-dimensional CW complex $X$ which is the union of its closed $n$-cells together with a continuous cellular map $f:X\rightarrow S_R$ whose restriction to each open cell is a homeomorphism. All tile types are $R$-complexes.

We now describe how to subdivide an $R$-complex $X$ with map $f:X\rightarrow S_R$, as described above. Recall that $R(S_R)$ is a subdivision of $S_R$. We simply pull back the cell structure on $R(S_R)$ to the cells of $X$ to create $R(X)$, a subdivision of $X$. This gives an induced map $f:R(X)\rightarrow R(S_R)$ that restricts to a homeomorphism on each open cell. This means that $R(X)$ is an $R$-complex with map $\phi_R \circ f:R(X)\rightarrow S_R$. We can iterate this process to define $R^n(X)$ by setting $R^0 (X) =X$ (with map $f:X\rightarrow S_R$) and $R^n(X)=R(R^{n-1}(X))$ (with map $\phi^n_R \circ f:R^n(X)\rightarrow S_R$) if $n\geq 1$.

We will use the term `subdivision rule' throughout to mean a finite subdivision rule of dimension $n$ for some $n$. As for traditional finite subdivision rules, we will describe an $n$-dimensional finite subdivision rule by a diagram giving the subdivision of every tile type, instead of by constructing an explicit complex.

\begin{figure}
\begin{center}
\scalebox{.5}{\includegraphics{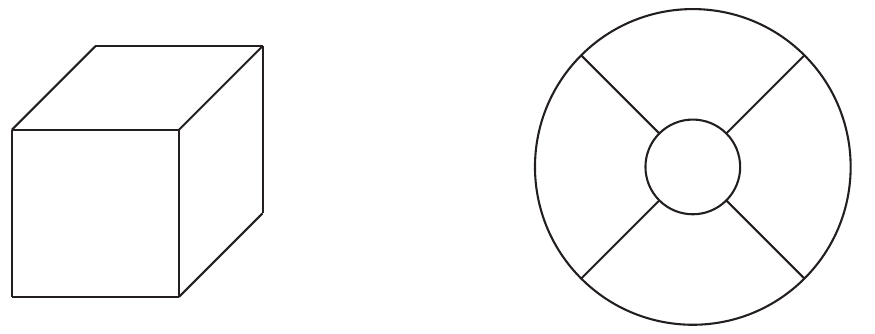}}
\caption{S(1)} \label{SOne}
\end{center}
\end{figure}

\begin{figure}
\begin{center}
\scalebox{.9}{\includegraphics{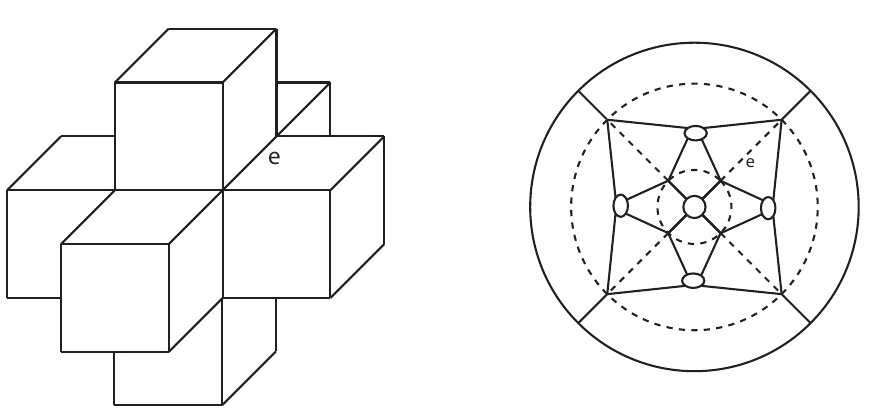}}
\caption{S(2)}\label{STwo}
\end{center}
\end{figure}

\begin{figure}
\begin{center}
\scalebox{.7}{\includegraphics{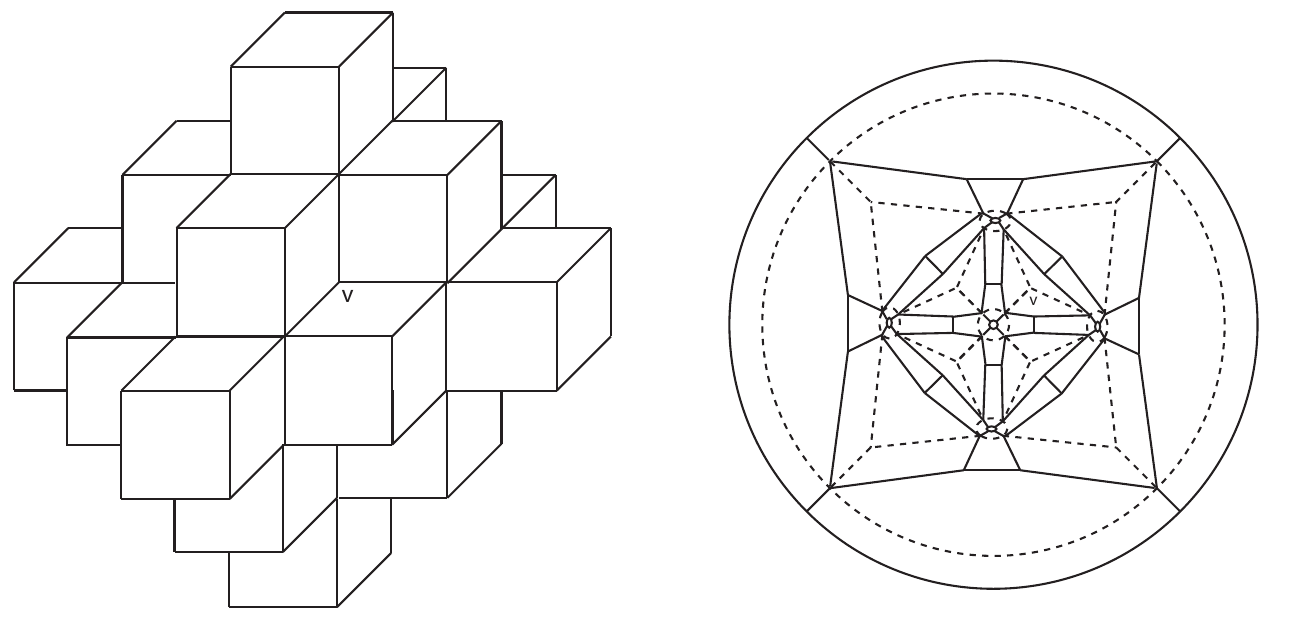}}
\caption{S(3)}\label{SThree}
\end{center}
\end{figure}

Our approach to finding subdivision rules in this paper and in others (see \cite{myself}, \cite{myself2}) is to take the boundary of balls in the universal cover. The universal cover of any manifold can be constructed recursively by taking a copy of the fundamental domain, gluing on fundamental domains to every exposed face of the original, and repeating. More specifically, let $B(0)$ be a single copy of the fundamental domain of an $n$-manifold $M$. Let $B(k)$ be the set of all fundamental domains that are distance $\leq k$ from $B(0)$ (in the word metric). Then for many groups and choices of generating sets, $S(k)=\partial B(k)$ will be a topological $(n-1)$-sphere for all $k$ or for $k$ sufficiently large. The cell structure from the fundamental domain gives a cell structure to $B(k)$ and thus to $S(k)$. This cell structure is a tiling. Thus, we get a sequence of tilings in which every tile or every group of tiles corresponds to an element of the fundamental group, and the entire group is represented at some point. We have drawn $S(1),S(2)$ and $S(3)$ for the 3-dimensional torus with the standard choice of generators in Figures \ref{SOne} to \ref{SThree}, shown in 3-space and also as a combinatorial tiling.

However, this sequence of tilings for a manifold is not necessarily created by a subdivision rule, because faces and edges are created and later covered up. To get a recursive structure, similar to hyperbolic 3-manifolds, we need to find a way to represent $S(k)$ (or a slightly modified version of it) as a subdivision of $S(k+1)$ (or a modified version of it).

\section{The $n$-torus}

We now show how to obtain a subdivision rule for the $n$-dimensional torus.

We will make the informal language of the introduction more rigorous. In the discussion that follows, let $I=[0,1]$, the unit interval. A $q$-cube is $I^q$, and a $p$-simplex is the convex hull of $p+1$ points in general position. Thus, a 1-cube is a line, a 2-cube is a square, and a 3-cube is a (standard) cube; a 1-simplex is a line, a 2-simplex is a triangle, and a 3-simplex is a tetrahedron, etc.

\begin{thm}\label{CubeTheorem}
The $n$-torus has a subdivision rule with $n$ tile types. The tile types are $p-1$-simplices cross $q$-cubes, where $1\leq p\leq n$ and $q=n-p$. Each such tile is subdivided into one $p-1$ simplex cross a $q$-cube and $2q$ $p$-simplices cross $q-1$ cubes.
\end{thm}

Before we begin the proof, Recall Figures \ref{SmallToriSubs1} and \ref{FourTorusSubs1} to see the tile types for $n=1,2,3$ and $4$.

\begin{figure}
\begin{center}
\scalebox{.8}{\includegraphics{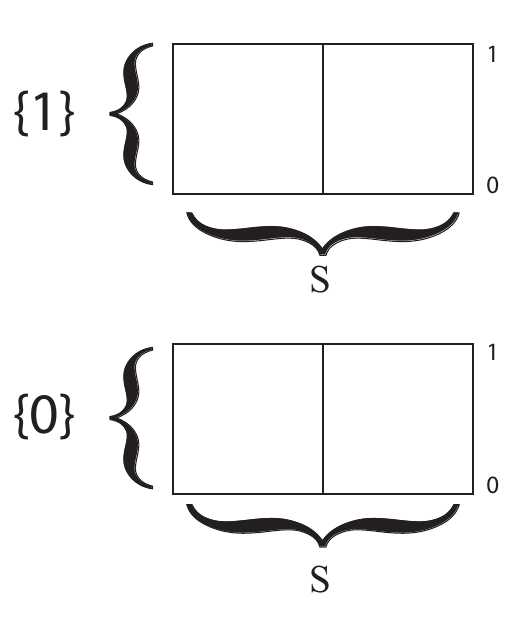}} \caption{The set $S\times \partial I^1\times I$, where $S$ is the union of two 1-simplices.} \label{PreCone}
\end{center}
\end{figure}

\begin{figure}
\begin{center}
\scalebox{.8}{\includegraphics{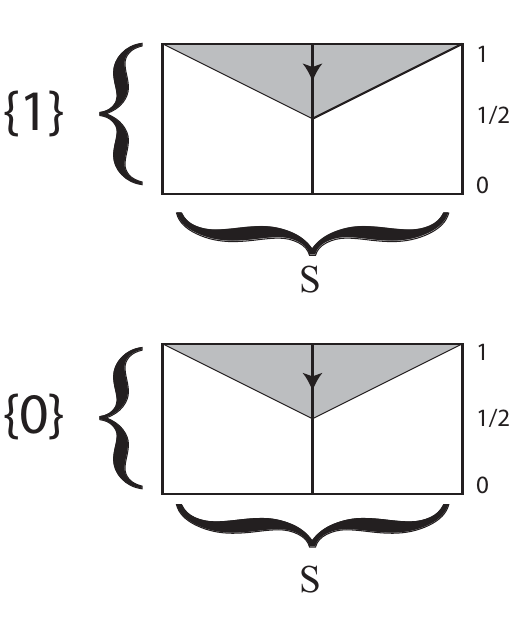}} \caption[An illustration of the homotopy.]{The same set as Figure \ref{PreCone}, after the homotopy. The grey
regions are the complement of the image of the homotopy; note that they are 2-simplices.} \label{Homotopy}
\end{center}
\end{figure}

\begin{figure}
\begin{center}
\scalebox{.8}{\includegraphics{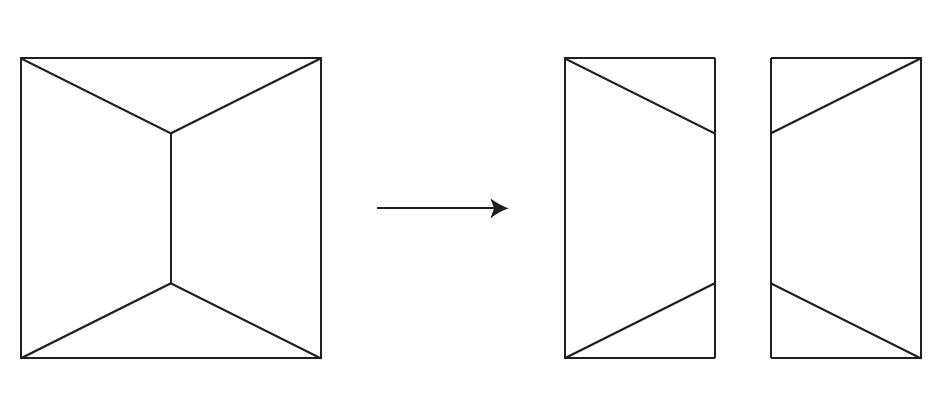}} \caption[An example of creating the subdivision rule.]{The same set as in Figure \ref{Homotopy}, after collapsing
$\partial I^1\times {0}$. Note that the original cell structure (2 squares) is contained in the new structure.} \label{FinalCone}
\end{center}
\end{figure}

\begin{figure}
\begin{center}
\scalebox{.8}{\includegraphics{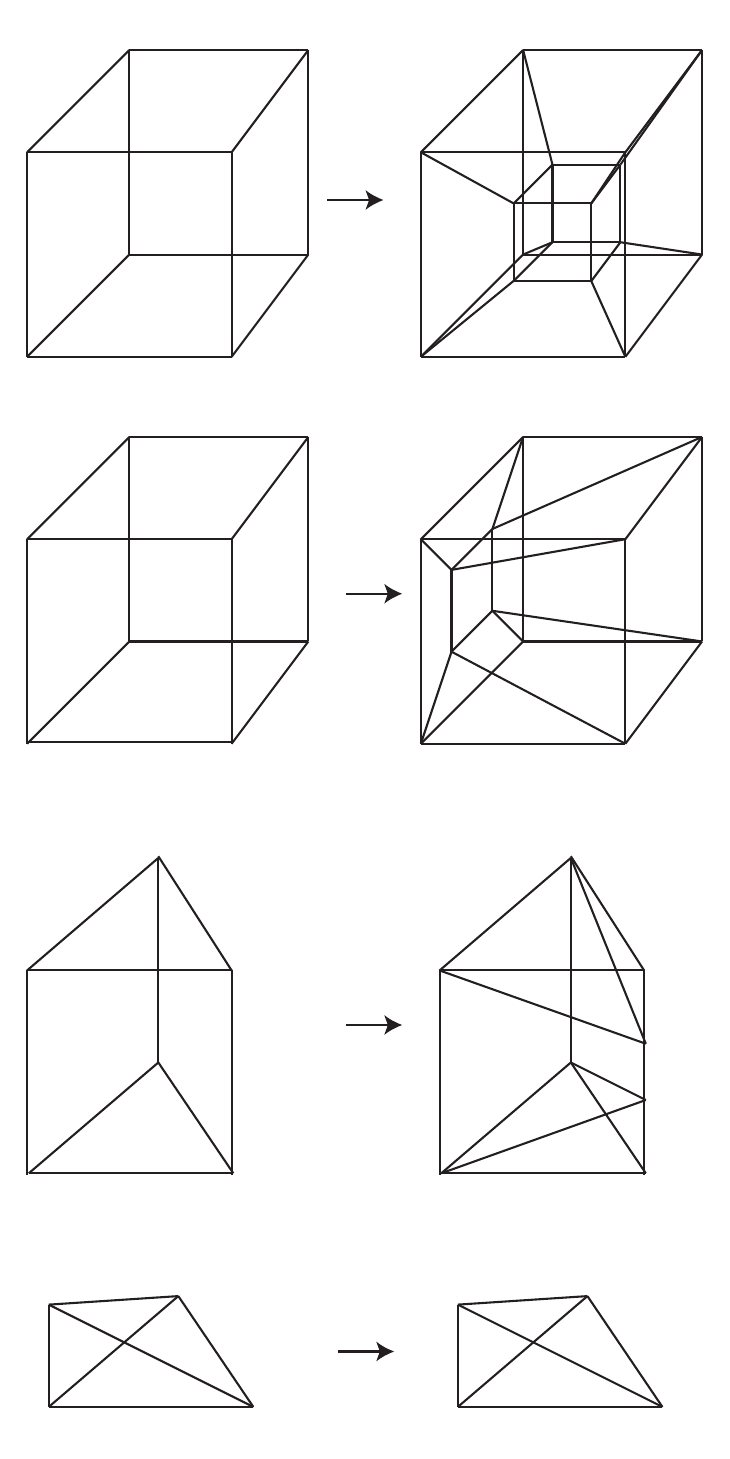}} \caption{The tile types for the four-torus.} \label{FourTorusSubs}
\end{center}
\end{figure}

\begin{proof}

The fundamental domain of the $n$-torus $\mathbb{T}^n=(S^1)^n$ is a hypercube of dimension $n$. If the generators of the fundamental group are
$y_1,...,y_n$, then every element of the fundamental group can be written uniquely as $y_1^{a_1} y_2^{a_2} ... y_n^{a_n}$.

Because our group is free abelian, the Cayley graph of the subgroup generated by any subset of the generators is contained in the Cayley graph of the fundamental group. Thus we can build the universal cover of these manifolds inductively from the universal covers of manifolds corresponding to subgroups.

We now describe how to explicitly construct the subdivision rule. It may help to follow along with the examples $n=1,2,3$ and 4 starting on page \pageref{Example}.

To construct the universal cover, we start with a single $n$-cube (i.e. $I^n$) and begin gluing on other cubes. Faces (or cells of codimension one) correspond to generators and inverses of generators. Assume an element represented by a cube is being glued on in some stage of creating the universal cover. Assume the element can be written as  $y_{k1}^{a_1}y_{k2}^{a_2} ... y_{kp}^{a_p}$, where this is a representation of minimal word length (so $1\leq k_1 \leq k_2 \leq...\leq k_p \leq n$ and $a_i \neq 0$). Then this element is contained in a subgroup of rank $p$. Let $q=n-p$. Then gluing on the cube corresponding to this element is accomplished by identifying some of its boundary to the previous stage of the universal cover. If we write the cube $I^n$ as $I^p \times I^q$, the boundary will be $\partial I^p \times I^q \cup I^p \times \partial I^q$.

Now, because the group element has $p$ geodesic paths into it (for instance, if $a_i$>0, going to  $y_{k1}^{a_1-1} y_{k2}^{a_2} ... y_{kp}^{a_p}$ and then going through the
$y_{k1}$-face to our element), our cube representing this element is glued onto $p$ faces in the previous stage at once. Each of the $p$ faces represents a generator, and if one generator is represented, its inverse is not, meaning that no pair of opposite faces is in the set of faces glued onto the universal cover.. The structure of the $n$-cube is such that every set of $p$ faces
containing no opposing pairs determines a unique $q$-cell which is common to all of them (so, for instance, in a 3-cube, three non-opposing faces intersect in a vertex, two in a line, and one in a square). If we project $I^n \subseteq \mathbb{R}^n$ down onto the subspace orthogonal to this cell, we see that this set of faces projects to the star of a vertex in $\partial I^p$. Call this star $S$. Note also that every vertex in the $p$-cube has an
opposite vertex, and the star of a vertex and its opposite have disjoint interiors and cover $\partial I^p$. Call the star of the opposite vertex $S^*$.

Thus, in gluing on $I^n$ via $\partial I^n$, we glue the boundary onto $A=S\times I^q$. The faces of the $\partial I^n$ that are not glued to anything can be written as $B=B_1\cup B_2= S^* \times I^q \cup I^p \times \partial I^q$. Recall that, to find a subdivision rule, we look at $S(k)$ (i.e. all exposed faces at stage $k$ of constructing a universal cover), and $S(k+1)$ (all exposed faces at stage $k+1$),and try to find the first as a subset of the second. Therefore, our goal is to find a cell structure for $A$ and $B$ such that $B$ is a refinement or subdivision of $A$. We use the standard simplicial decomposition of the $p$-cube (found, for instance, in \cite{Rudin}, Exercise 10.18), which we now describe.

$I^p$ is covered by the $p!$ simplices $\{[0,e_{\sigma(1)},e_{\sigma(1)}+e_{\sigma
(2)},e_{\sigma(1)}+e_{\sigma (2)}+...+e_{\sigma(p)}]| \sigma\in \Sigma_p \}$, each of which has disjoint interior. Here, $e_i$ is the unit vector in the $i$-th
direction. The symbol $[p_0,p_1,...,p_k]$ is defined to be $\tau(Q^k)$, where $\tau$ is the affine map $\tau(x_1,...,x_k)=p_0+\Sigma x_i(p_i-p_0)$, and $Q^k$ is
the standard simplex $\{(x_1,...,x_k)|0\leq x_i$ for all i$, x_1+...+x_n \leq 1\}$. Each of these simplices has sub-simplices defined by deleting intermediate terms (so $[0]\subseteq [0,e_1] \subseteq [0,e_1,e_2]$, for instance).

Recall that switching two terms in the simplex (i.e. changing $[p_0, p_1, p_2]$ to $[p_2,p_1,p_0]$) gives a different map from $Q^k$ with opposite orientation but the same image as the original map. If $\tau_1$ and $\tau_2$ are the map corresponding to the original simplex and the `flipped' simplex, then $\tau_2^{-1}\tau_1$ is an orientation-reversing simplicial map.

We use this to define an involution on our simplicial cube. Define this map by switching $0$ and $e_{\sigma(1)}+e_{\sigma (2)}+...+e_{\sigma(p)}$ in every simplex. This is a simplicial map that is the identity on all subsimplices not containing 0 or $e_{\sigma(1)}+e_{\sigma (2)}+...+e_{\sigma(p)}$. Any subsimplex that contains one of those points is sent to an opposing subsimplex that contains the other point. The existence of this map shows, in particular, that the set of all closed simplices in $\partial I^p$ containing 0 is simplicially isomorphic to the set of all closed simplices in $\partial I^p$ containing $e_{\sigma(1)}+e_{\sigma(2)}+...+e_{\sigma(n)}$. This means that if $S$ and $S^*$ are given the simplicial structure they inherit from $I^p$, they are isomorphic.

This means that $A$ and $B_1$ have the same cell structure. If $q=0$, $B_1=B$ and $A$ and $B$ have the same cell structure, so our subdivision rule can be the
identity.

If $q\neq 0$, it's slightly more difficult. We still give $A$ and $B_1$ the simplicialized structure explained above, and give $B_2$ the structure of $I^p \times \partial I^q$, where $I^p$ is given the simplicial structure of earlier. We show that $B$ contains $A$ as a subcomplex, with $\partial A \subseteq \partial B$. In the discussion that follows, it will be helpful to follow along with Figures \ref{PreCone}-\ref{FinalCone} for the case $p=2$, $q=1$. Figure \ref{FourTorusSubs} gives more examples with less explanation.

So, pick a $(p-1)$-simplex $\Delta^{p-1}$ in $S \subseteq \partial I^p$. If we consider the center vertex of $S$ as $0$ and the center vertex of $S^*$ as
$e_{\sigma(1)}+e_{\sigma (2)}+...+e_{\sigma(p)}$, then there is a unique $p$-simplex $\Delta^p$ defined by adjoining $e_{\sigma(1)}+e_{\sigma
(2)}+...+e_{\sigma(p)}$ to $\Delta^{p-1}$. So, our goal is to show that $\Delta^{p-1}\times I^q\subseteq \Delta^{p-1}\times I^q \cup \Delta^p \times \partial
I^q$. Patching together these simplices will show that $A\subseteq B_1\cup B_2=B$.

To do so, note that $I^q$ is just the cone over the boundary of $I^q$ (as a set, not as a complex). Thus, we look at $\Delta^{p-1} \times \partial I^q \times I$,
which we will eventually collapse. Each face in $\partial I^q$ is a $(q-1)$-cube. Given a specific face, we can embed the product $\Delta^{p-1}\times I^{q-1}
\times I$ in $\mathbb{R}^{p+q-1}$ as $\{(x_1,...,x_{p-1},y_1,...,y_{q-1},z)|0\leq x_i$ for $1\leq i \leq p-1, 0\leq y_j \leq 1$ for $1\leq j \leq q-1,0 \leq z \leq 1 x_1+x_2+...+x_{p-1}\}$. Call this set $C$.

Define a family of maps $f_t:C \rightarrow C$ by $$f_t(x_1,...,x_{p-1},y_1,...,y_{q-1},z)=(x_1,...,x_{p-1},y_1,...,y_{q-1},z(1+(x_1+...+x_{p-1}-1)\frac{t}{2})).$$ This defines an invertible homotopy
(basically dragging down the corner of the top copy of the simplex along the $z$-axis).

Note that $$f_1(C)=\{(x_1,...,x_{p-1},y_1,...,y_{q-1},z)|z \leq \frac{1}{2}+\frac{x_1+...+x_{p-1}}{2}\}$$ and this is the same as
$$\{(x_1,...,x_{p-1},y_1,...,y_{q-1},z)|x_1+...+x_{p-1}+(2-2z)\geq 1\}.$$ The closure of its complement in $C$ is $I^{q-1}$ cross a $p$-simplex defined by
$x_1+...+x_{p-1}+2(1-z) \leq 1$, where $0 \leq x_i, z \leq 1$. Thus we can write $C$ as the union of $f_1(C)\cong C$ and $C\setminus f_1(C)\cong \Delta^p \times
I^{q-1}$. The boundary of $f_1(C)\cup(C\setminus f_1(C))$ is clearly the same as $C$, just with a more complex cell structure, i.e. a subdivision. If we now
collapse to get the cone structure mentioned earlier, the simplex we just obtained is not affected, and we still have a subdivision. Patching together all faces
in $\partial I^q$ shows that $\Delta^{p-1}\times I^q\subseteq \Delta^{p-1}\times I^q \cup \Delta^p \times \partial I^q$, as desired; since the homotopy fixed all $y$-coordinates, the subdivisions on each face of $\partial I^q$ match up. Finally, we glue together all the simplices to show that $A\subseteq B_1 \cup B_2$. Note that the center vertex of $S$ was the corner of each simplex sent to the origin in our embedding above, so when we glue our simplices together, all of those vertices are identified, and we have a well-defined subdivision rule.

There is one problem that may have arisen. First of all we are assuming that $A$ is formed of $(p-1)$-simplices crossed with $q$-cubes. We also have a structure
for $B$; but in the next stage of subdivision (or of constructing the universal cover), the new $A$'s are formed from the old $B$'s. Do they have the right
structure? Well, note that $B_1$ was given the same structure as $A$, and the subtiles of $B_1$ correspond to those elements of $Z^n$ that stay in the same
subgroup of rank $p$. $B_2$ represents those elements that land in a subgroup of order $p+1$, and these are given the structure of $n$-cubes split into
$p$-simplices cross $(q-1)$ cubes, and the $p$-simplices are grouped about the correct vertex. So, the cell structure is consistent.

\end{proof}

We now find the subdivision rules explicitly for $n=1,2,3$ and $4$. \label{Example}

\begin{figure}
\begin{center}
\scalebox{1.4}{\includegraphics{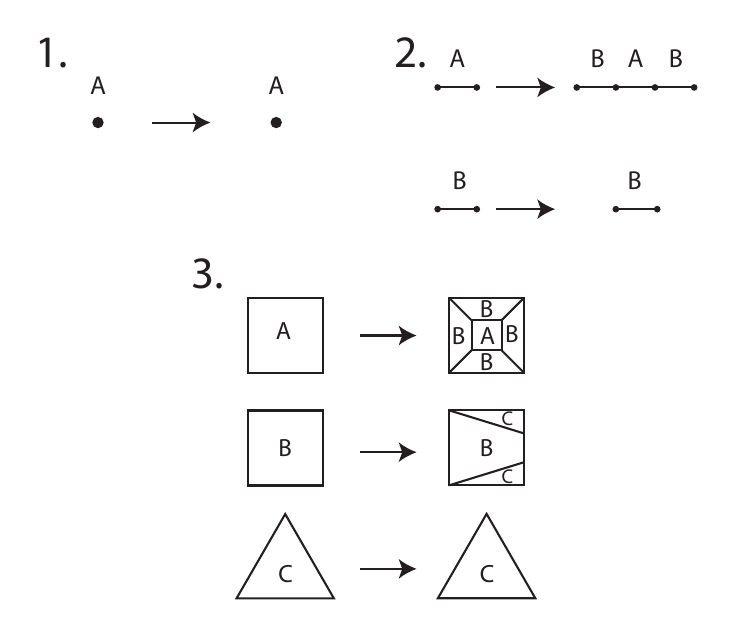}} \caption[The subdivision rules for the 1-,2- and 3-torus.]{1. The subdivision rule for the 1-torus, i.e. the
circle. 2. The subdivision rule for the 2-torus, i.e. the standard torus. 3. The subdivision rule for the 3-torus.} \label{SmallToriSubs}
\end{center}
\end{figure}

For $n=1$, the universal cover is a line, $B(n)$ is $2n+1$ line segments, and $S(n)$ is two points. The subdivision rule is shown in Figure \ref{SmallToriSubs}.

Note that the only tile type is a point (i.e. 0-simplex cross a 0-cube), which is subdivided into one 0-simplex cross a 0-cube.

For $n=2$, the fundamental domain is a square, $B(n+1)$ is a topological disk, and $S(n)$ is a topological circle. The subdivision rule is shown in Figure
\ref{SmallToriSubs}.

Type A is a line (i.e. a 0-simplex cross a 1-cube). It's subdivided into one line( a 0-simplex cross a
1-cube) and 2 more lines (1-simplices cross a 0-cube), just as the formula predicts.

Type B is a line(a 1-simplex cross a 0-cube), and represents half of a group element. Two B tiles form the star of a vertex in the boundary of the 2-cube (a square), and
the subdivision rule for type B is the identity, just as the formula predicts.

\begin{figure}
\begin{center}
\scalebox{.5}{\includegraphics{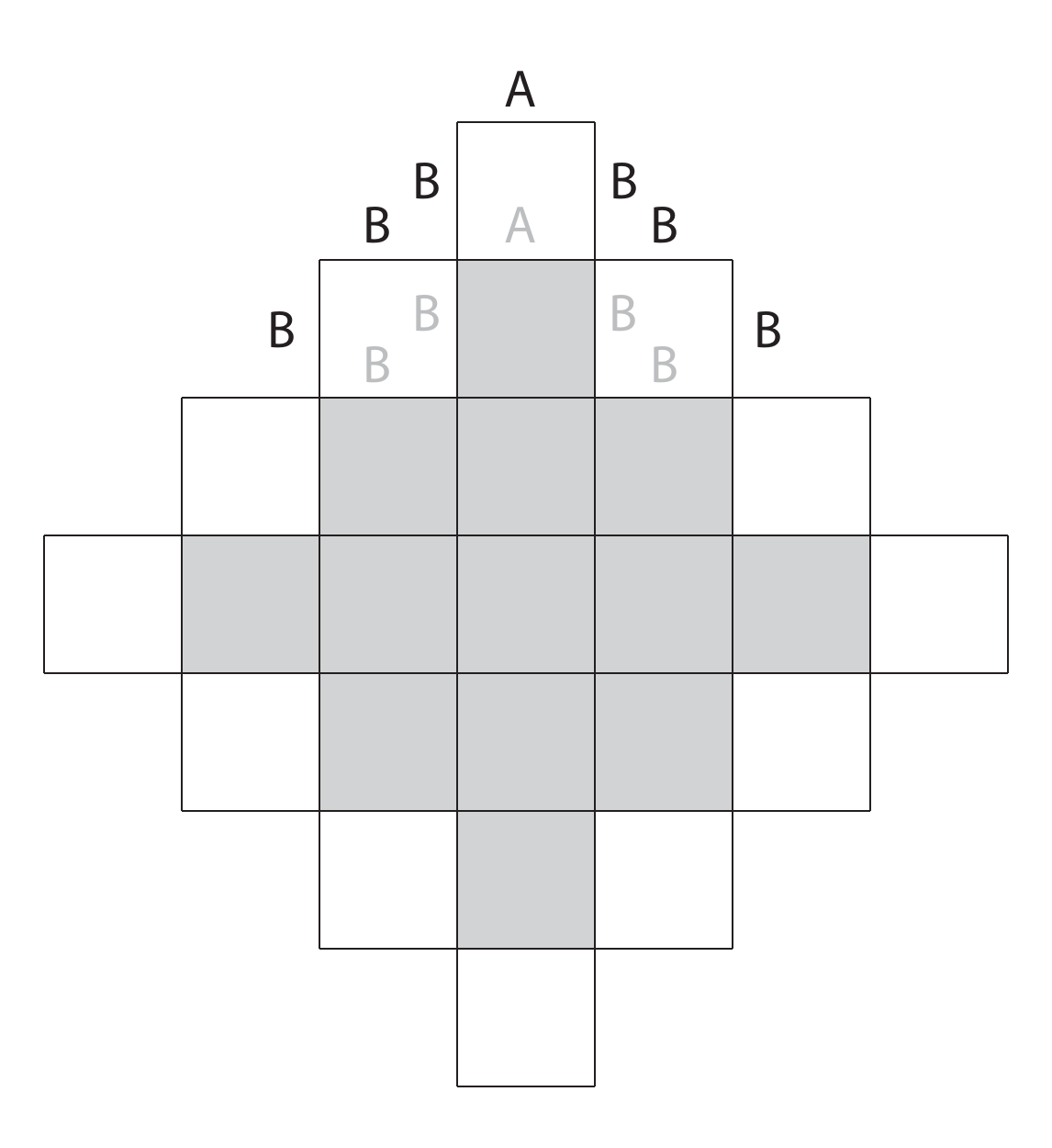}} \caption[The universal cover of the 2-torus.]{White tiles correspond to $B(3)$, and grey tiles to $B(2)$.
Note that A tiles correspond to elements furthest from the origin.} \label{TorusUniversal}
\end{center}
\end{figure}

As you can see in Figure \ref{TorusUniversal}, A tiles correspond to the four `ends' of $B(n)$, or, the group elements contained in a subgroup generated by
exactly one of the standard generators, while B tiles correspond to elements that must be written using both generators. Notice how two neighboring $B$ tiles form a corner that is covered up by one square fundamental domain.

For $n=3$, we have the 3-torus, whose universal cover is shown being constructed in Figures \ref{SOne} to \ref{SThree}. Notice in these figures that new cubes are glued onto a single face, two neighboring faces, or three faces forming a corner. These correspond to elements whose minimal word-length representations use one, two, or three generators, respectively. Type A (corresponding to a single face) is a square (or 0-simplex cross a 2-cube), and is subdivided into one square (a 0-simplex cross a 2-cube) and 4 other squares (or 1-simplices cross 1-cubes).

Type B is a square (thought of as a 1-simplex cross a 1-cube). It is subdivided into a square (a 1-simplex cross a 1-cube), and 2 triangles (or 2-simplices cross a 0-cube). Two type B tiles correspond to the star of a vertex in the boundary of the 2-cube, which is then crossed by $I$.

Note that this tile shows us what happens with the homotopy portion of Theorem \ref{CubeTheorem}. We start with $S$, the star of a vertex in the boundary of
$I^p=I^1$ with a simplicial structure (namely, two edges of a square sharing a vertex, each edge considered as a 1-simplex). We cross this with $I^q$=$I^1$ to get two square
sharing a face. This is $A$. Then, we write $A$ as a quotient of $S\times \partial I^1\times I$, so we get Figure \ref{PreCone}.

On each component (one corresponding to $S\times{0}\times I$, one corresponding to $S\times{1}\times I$), we pull down the center line by our homotopy to get
Figure \ref{Homotopy}. Collapsing $\partial I^1\times \{0\}=\{0,1\}\times\{0\}$ to a point to get a cone, we have the image in Figure \ref{FinalCone}.

Type C is a triangle (a 2-simplex cross a 0-cube), which is subdivided by the identity. This corresponds to a cube glued onto three adjoining faces. Each tile of type C represents one sixth of a group element, and six together form the star of a vertex in $\partial I^3$ when it is given a simplicial structure.

\begin{figure}
\begin{center}
\scalebox{.7}{\includegraphics{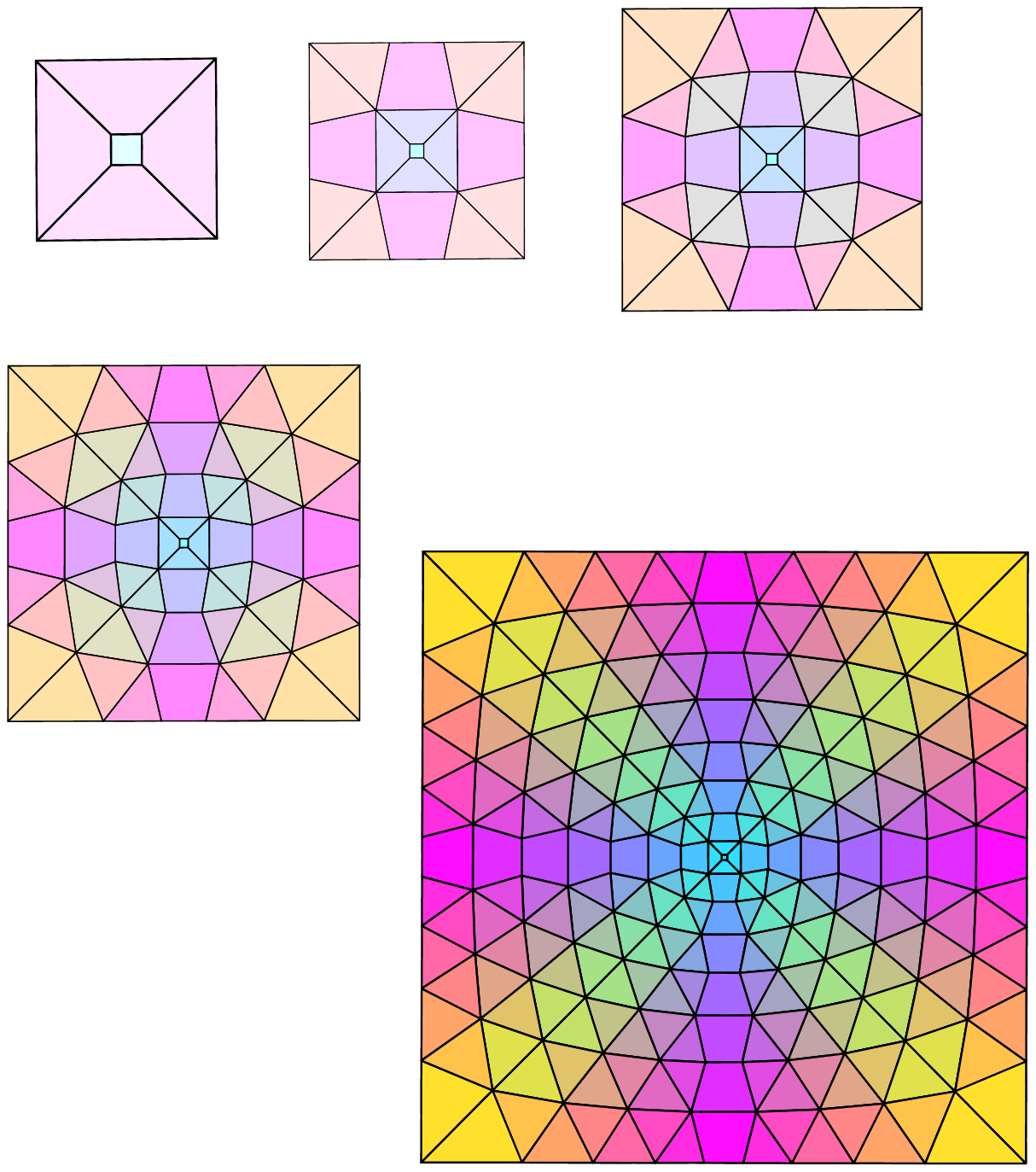}} \caption{Several subdivisions of a type A tile for the 3-torus.} \label{CircleTorus}
\end{center}
\end{figure}

Several subdivisions of an A tile are shown in Figure \ref{CircleTorus}. This picture was created with Ken Stephenson's Circlepack \cite{Circlepak}. The pictures are only combinatorial subdivisions of each other; they can't be overlaid with vertices matching up. This is because the subdivision rule is not conformal. For more on the connection between conformality and circle packings, see \cite{French}.

Finally, the tile types for the four-torus are shown in Figure \ref{FourTorusSubs}.

Note the homotopies in each of these tiles. In the first tile type, the homotopy is shrinking the cube in the center and dragging the cell structure with it. In the second, we only shrink down a 2-dimensional face, again dragging everything along with it. In the third tile type, we shrink an edge, dragging along the cell structure with it. Finally, in the fourth tile type, there is nothing to drag.

\section{Future Work}

We hope to find many more subdivision rules for higher-dimensional manifolds, especially hyperbolic manifolds. Cannon and Swenson have shown \cite{hyperbolic} that hyperbolic $n$-manifolds have some sort of subdivision rule. We hope to find more explicit examples.

Also, as mentioned in the introduction, Hersonsky and others have studied extremal length for three dimensional tilings (see \cite{saar}, \cite{Spherepackings}, and \cite{Spherelackings}). Do 3-dimensional subdivision rules for hyperbolic 4-manifolds satisfy a condition on extremal length similar to conformal 2-dimensional subdivision rules?

\bibliographystyle{plain}
\bibliography{CubePaper}

\end{document}